\documentclass{amsproc}

\usepackage{amssymb}

\newtheorem{theorem}{Theorem}[section]

\newtheorem{definition}[theorem]{Definition}

\newtheorem{remark}[theorem]{Remark}

\newcommand{\coleq}{{\colon}{=}\;}
\newcommand{\er}{\mathbb R}

\newcommand{\en}{\mathbb N}
\newcommand{\codim}{\mathrm{codim}}
\newcommand{\hj}{J^r(M,N)}
\newcommand{\hjs}{J^r_x(M,N)}
\newcommand{\hjt}{J^r(M,N)_{\bar x}}
\newcommand{\hjst}{J^r_x(M,N)_{\bar x}}

\newcommand{\nj}{\tilde J^r(M,N)}
\newcommand{\qj}{Q\! J^r(M,N)}

\newcommand{\qjst}{Q\! J^r_x(M,N)_{\bar x}}

\newcommand{\iter}{\underbrace{T\dots T}_{\textrm{\tiny $r$--times}}}
\newcommand{\iterm}{\underbrace{T\dots T}_{\textrm{\tiny $r-1$--times}}}
\newcommand{\itea}{\underbrace{T\dots T}_{\textrm{\tiny $a$--times}}}
\newcommand{\ites}{\underbrace{T\dots T}_{\textrm{\tiny $s$--times}}}
\newcommand{\itebm}{\underbrace{T\dots T}_{\textrm{\tiny $b$--times}}M}
\newcommand{\itebmm}{\underbrace{\scriptstyle T\dots T }_{\textrm{\tiny $b$--times}}M}
\newcommand{\ps}{\pi^s}
\newcommand{\psb}{\pi^s_b}
\newcommand{\aps}{{}_a^{}\!\pi^s_{}}
\newcommand{\apsb}{{}_a^{}\!\pi^s_b}
\newcommand{\apjb}{{}_a^{}\!\pi^1_b}

\newcommand{\gel}{\mathrm{GL}}
\newcommand{\dej}{\textrm{d}}

\begin{document}


\author{Miroslav Kure\v{s}}
\title{On some directions in the development of jet calculus}

\address
{Miroslav Kure\v{s}, Institute of Mathematics, Brno University of Technology,
Technick\'a 2, 61669 Brno, Czech Republic}
\email
{kures@fme.vutbr.cz}

\maketitle

\begin{abstract}
Two significant directions in the development of jet calculus are showed.
First, jets are generalized to so-called quasijets.
Second, jets of foliated and multifoliate manifold morphisms are presented.
Although the paper has mainly a survey character, it also includes new results: jets modulo multifoliations are introduced
and its relation to $(R,S,Q)$-jets is demonstrated. 
\end{abstract}

\section{Introduction} 
In this note, we present directions in a generalization of a notion of the jet.
Thus, roughly speaking, the paper can be read as observations about a development of jet calculus.
We start by nonholonomic and semiholonomic jets, which were defined already by Ehresmann; as a more general concept,
quasijets were introduced and studied at first by Pradines in \cite{PRA}.
So, we recall some facts about nonholonomic jets
and quasijets in Section~1. In Section~2, we mention fundamentals from the theory of foliations and
show various multifoliate structures.
Further, we initialize jet formalism for such multifoliate structures.
Here, we have two inspirations:
Ikegami's paper \cite{IKE} about jets modulo foliations and the concept of $(R,S,Q)$-jet,
see e.g. \cite{KMS} or \cite{DOK}. We present a way to a generalization and
unification (in a way) of both concepts as our new results.
Section 3 is devoted to interactions between mentioned generalizations, to a Weil algebras approach and a concept of weighted jets.
All manifolds and maps are assumed to be of class $C^\infty$.

\section{From holonomic jets via nonholonomic jets to quasijets}

\subsection{Holonomic jets}

Jets (holonomic jets) are commonly known as certain equivalence classes
of smooth maps between manifolds, which are represented by Taylor polynomials. First, we precise the definition. 
Let $M$ and $N$ be two manifolds. Then
two maps 
$$
f\colon M\to N,\quad g\colon M\to N
$$
are said to {\it determine the same
$r$--jet} at $x\in M$, if for every curve 
$$
\gamma\colon\mathbb R \to M
\;\;\textrm{with}\;\;
\gamma(0)=x
$$
the curves 
$$
f\circ\gamma
\;\;\textrm{and}\;\;
g\circ\gamma
$$
have the $r$--th
order contact at $0\in\er$. In such a case we write 
$$
j^r_xf=j^r_xg
$$
and an equivalence class of this relation is called an {\it $r$--jet} of $M$ into
$N$. The set of all $r$--jets of $M$ into $N$ is denote by $\hj$. If the
{\it source} of a $r$--jet is $x\in M$ and the {\it target} of this jet is
$\bar x=f(x)\in N$,
then
$$
\alpha\colon j^r_xf\mapsto x
\;\;\textrm{and}\;\;
\beta\colon j^r_xf\mapsto \bar x
$$
are projections of
fibered manifolds 
$$
\alpha\colon\hj\to M, \qquad \beta\colon\hj\to N.
$$
Further,
by 
$$
\hjs
\;\;\textrm{or}\;\;
\hjt
$$ we mean the set of all $r$--jets of $M$ into $N$ with the
source $x\in M$ or with the target $\bar x\in N$, respectively, and we write
$$
\hjst=\hjs\cap\hjt.
$$
As $r$--th order contact of maps is preserved under
composition, we define the composition of $r$--jets as the $r$--jet of
composed map.

Let $p\colon Y\to M$ be a fibered manifold. The set $J^rY$ of all $r$--jets
of the local sections of $Y$ is called the {\it $r$--th jet prolongation of
$Y$} and $J^rY\subset J^r(M,Y)$ is a closed submanifold. (If $Y\to M$ is a
vector bundle, then $J^r Y\to M$ is also a vector bundle.) Let $q\colon Z\to N$
be another fibered manifold and $f\colon Y\to Z$ a fibered bundle morphism with
the property that the base map $f_0\colon M\to N$ is a local diffeomorphism.
There is an induced map 
$$
J^r(f_0,f)(X)\coleq
j^r_{\beta(X)}f\circ X\circ j^r_{f_0(\alpha(X))}f_0^{-1}
$$
for $X\in J^r(M,Y)$. If we restrict it to local sections, we obtain a map
denoted by 
$$
J^rf\colon J^rY\to J^rZ
$$
which is called the {\it $r$--th jet
prolongation of $f$}.

It is clear that the trivial choice $Y=M\times N$ 
in the jet prolongation yields jets of mappings from $M$ to $N$.

\subsection{Nonholonomic jets}

For $r=1$, the set of nonholonomic 1--jets 
$$
\tilde J^1(M,N)\coleq J^1(M,N).
$$
By induction, let
$\alpha\colon\tilde J^{r-1}(M,N)\to M$ denote the source projection and
$\beta\colon\tilde J^{r-1}(M,N)\to N$ the target projection of
$(r-1)$--th nonholonomic jets. Then $X$ is said to be a {\it nonholonomic
$r$--jet} with the source $x\in M$ and the target $\bar x\in N$, if there
is a local section 
$$
\sigma\colon M\to\tilde J^{r-1}(M,N)
$$
such that
$$
X=j^1_x\sigma
\;\;\textrm{and}\;\; 
\beta(\sigma(x))=\bar x.
$$ 
There is a natural
embedding $\hj\subset\nj$.

Every $X\in\nj$ induces a map
$$
\mu X\colon (\iter M)_x \to (\iter N)_{\bar x}
$$ 
in the following way.
For $r=1$ and $X=j^1_xf$ is $\mu X$ defined as $T_xf$. By induction,
let $X=j^1_x\sigma$ for a local $\alpha$--section
$$
\sigma\colon M\to \tilde J^{r-1}(M,N).
$$ 
Then $\sigma(u)\in J_u^{r-1}(M,N)$,
$$
\mu(\sigma(u))\colon(\iterm M)_u\to(\iterm N)_{\beta(\sigma(u))}
$$
and we
put 
$$
\mu X=T_x\mu(\sigma(u)).
$$
The constructed map
$$
\mu X\colon (\iter M)_x \to (\iter N)_{\bar x}
$$
is a vector bundle
morphism with respect to all vector bundle structures
$\iter\to\iterm$. However, $\mu X$ is not an entirely general
vector bundle morphism of this type.

Let $p\colon Y\to M$ be a fibered manifold. We construct the {\it $r$--th jet
nonholonomic prolongation of $Y$} denoted by $\tilde J^rY$ as the set of all
nonholonomic $r$--jets of the local sections of $Y$. The construction of the
{\it $r$--th jet nonholonomic prolongation of $f$} for a fibered bundle
morphism $f\colon Y\to Z$ with the property that the base map
$f_0\colon M\to N$ is a local diffeomorphism is analogous to the holonomic
case.

\subsection{Quasijets}

We introduce the following denotation of projections in the iterated tangent
bundle $\iter M$. For every $s$, $0 < s\le r$, we denote by
$$
\ps\colon \ites M\to M
$$ the canonical projection to the base.
Further, we denote
$$
\psb\coleq \ps_{\itebmm}
\colon
\tilde T^s(\itebm)\to \itebm
$$
projection with $\itebm$ as the base space,
$$
\aps\coleq \itea \ps
\colon
\itea (\ites M)\to \itea M
$$
induced projection originating by the posterior application of the
functor $\itea$, 
$$
\apsb\coleq \itea\ps_{\itebmm}
$$
the general case
containing both previous cases.
If $a$ or $b$ equal zero, we do not write them.

Let $x\in M$, $\bar x\in N$. A map
$$
\phi\colon (\iter M)_x \to (\iter N)_{\bar x}
$$
is said
to be a {\it quasijet} of order $r$ with the source $x$ and the
target $\bar x$, if it is a vector bundle morphism with respect to all
vector bundle structures 
$$
\apjb\colon (\iter M)_x \to (\iterm M)_x
$$
and 
$$
\apjb\colon (\iter N)_{\bar x} \to (\iterm N)_{\bar x},
$$
$a+b=r-1$. The set of all such
quasijets is denoted by $\qjst$ and $\qj$ means the set of all quasijets
from $M$ to $N$.

There is a bundle structure $\qj\to M\times N$ and, analogously to $J^r$,
the set $Q\! J^rY$ of all $r$--jets of the local sections of
a fibered manifold $Y\to M$ $Y$ is called the {\it $r$--th quasijet
prolongation of $Y$}.
We compose quasijets as maps.
Further, let $q\colon Z\to N$
be another fibered manifold and $f\colon Y\to Z$ a fibered bundle morphism
with the property that the base map $f_0\colon M\to N$ is a local
diffeomorphism.
There is an induced map 
$$
Q\! J^r(f_0,f)(X)=
j^r_{\beta(X)}f\circ X\circ j^r_{f_0(\alpha(X))}f_0^{-1}
$$
for $X\in Q\! J^r(M,Y)$.
The composition denoted the composition of quasijets, where
the holonomic jets
$j^r_{\beta(X)}f$, $j^r_{f_0(\alpha(X))}f_0^{-1}$
are considered as quasijets by the use of
the map $\mu$ from Section~2.
If we confine ourselves to local sections, we obtain a map
denoted by $Q\! J^rf\colon Q\! J^rY\to Q\! J^rZ$ which is called
the {\it $r$--th quasijet prolongation of $f$}.

\section{Jets preserving foliated structures}

\subsection{Fundamentals from the theory of foliations}

Let $M$ be a $m$-dimensional smooth manifold, $m=p+q$, $m\in\en$, $p,q\in\mathbb N\cup\{ 0\}$, $(x,y)=(x^1,\dots,x^p,y^1,$ $\dots,y^q)\in\er^p\times\er^q=\er^m$. For constants $\bar c\in\er^p$, $c\in\er^q$,
we consider spaces $\er^q_{\bar c}=\{(x,y)\in\er^m; x^1=\bar c^1,\dots,x^p=\bar c^p\}$ and $\er^p_c=\{(x,y)\in\er^m; y^1=c^1,\dots,y^q=c^q\}$. Intersections of $\er^q_{\bar c}$ and $\er^p_c$ with open sets (balls) with respect to the standard topology are denoted by $P^q_{\bar c}$ and $P^p_c$ and called the {\it $(\bar c,q)$-coplaque} and the {\it $(p,c)$-plaque} in $\er^m$. Suppose that $\mathcal F=\{L_t\}_{t\in J}$ is a partition of $M$ into connected subsets, $M=\bigcup_{t\in J} L_t$, $L_t\cap L_s=\emptyset$ for $t\ne s$. Further, we consider a {\it foliated atlas} on $M$, i.e., a collection $\{U_i,\varphi_i\}_{i\in I}$, $\varphi_i=\alpha_i\times\beta_i$, $\alpha_i\colon U_i\to\er^p$, $\beta_i\colon U_i\to\er^q$, of charts
satisfying
\begin{enumerate}
\item[(i)] $\{U_i\}_{i\in I}$ is a cover of $M$ by open sets
\item[(ii)] each connected component of $L_t\cap U_i$ (for all $i\in I$, $t\in J$) is mapped by $\varphi_i$ onto
an $(p,c)$-plaque in $\er^m$, i.e., for $u\in U_i$
\begin{eqnarray*}
x^a=\alpha^a_i(u) &&\qquad a=1,\dots,p\\
y^b=\beta^b_i(u)=c^b &&\qquad b=1,\dots,q
\end{eqnarray*}
\item[(iii)] transition functions $\varphi_{ij}=\varphi_j\circ\varphi_i^{-1}$ on $U_i\cap U_j$, $\varphi_{ij}=\alpha_{ij}\times\beta_{ij}$, send $(p,c)$-plaques onto $(p,c)$-plaques, i.e.
\begin{eqnarray*}
x^a=\alpha^a_{ij}(x,y) &&\qquad a=1,\dots,p\\
y^b=\beta^b_{ij}(y) &&\qquad b=1,\dots,q.
\end{eqnarray*}
\end{enumerate}
Then $\mathcal F$ is called the {\it foliation} of $M$ of {\it dimension} $p$ and {\it codimension} $q$,
$L_t$, $t\in J$ {\it leaves} of $\mathcal F$ and $M$ the {\it foliated manifold} written shortly by $(M,\mathcal F)$.
Trivial cases arise for $p=0$, $q=m$ (leaves = points) and for $p=m$, $q=0$ (the unique leaf = $M$).

Let $\mathcal F$, $\mathcal F^\prime$ be two foliations of $M$ with dimensions $p$ and $p^\prime$. Then $\mathcal F^\prime$ is called a {\it subfoliation} of $\mathcal F$ and $\mathcal F$ is called a {\it superfoliation} of $\mathcal F^\prime$, denoted by
$\mathcal F^\prime\preceq \mathcal F$, if the following conditions hold:
\begin{enumerate}
\item[(i)] $0\le p^\prime\le p\le m$ 
\item[(ii)] for any leaf $L^\prime$ of $\mathcal F^\prime$, there exists a leaf $L$ of $\mathcal F$ such that $L^\prime\subseteq L$, and the restriction of $\mathcal F^\prime$ on a leaf $L$ of $\mathcal F$ is a foliation of dimension $p-p^\prime$ of $L$.
\end{enumerate}
The relation $\preceq$ is an order in the set of foliations of $M$.

Fibered manifolds are canonically foliated, their fibers can be viewed as leaves. On the other hand, there exist manifolds, which are foliated but not fibered.

\subsection
{Transversality of maps, transversality of foliations and induced multifoliations}

Let $\Delta$ be an integer greater than 1. Let us consider manifolds $H_\delta$, $\delta=1,\dots,\Delta$, and $M$.
Let $f_\delta\colon H_\delta\to M$, $\delta=1,\dots,\Delta$, be (smooth) maps.

We take an arbitrary non-empty subset $E\subseteq\{1,\dots,\Delta\}$ and denote by $\mathrm{Im}f_E$ the intersection of all images of $f_\epsilon$, $\epsilon\in E$.

For $u_E\in\mathrm{Im}f_E$ and every $\epsilon\in E$,
let $(Tf_\epsilon)_{u_E}$ denote the image of the tangent map to $f_\epsilon$ in $u_E$; 
tangent vectors belonging to $(Tf_\epsilon)_{u_E}$ generate a vector subspace
of $T_{u_E}M$; we denote it by $\langle(Tf_\epsilon)_{u_E}\rangle$. Further, we denote by
$\langle \bigcup\limits_E (Tf_\epsilon)_{u_E} \rangle$ 
the vector space generated by the union of vectors in all $(Tf_\epsilon)_{u_E}$, $\epsilon\in E$, and
by $\langle \bigcap\limits_E (Tf_\epsilon)_{u_E} \rangle$ the vector space generated by vectors belonging to the intersection of all $(Tf_\epsilon)_{u_E}$, $\epsilon\in E$.

For simplicity, we consider only maps for which vector spaces above
have constant dimensions  for all $u_E\in\mathrm{Im}f_E$.

Now, it is evident that for every chosen $\epsilon_0\in E$
$$
0\le \dim\langle \bigcap\limits_E (Tf_\epsilon)_{u_E} \rangle
\le \dim\langle (Tf_{\epsilon_0})_{u_E} \rangle
\le \dim\langle \bigcup\limits_E (Tf_\epsilon)_{u_E} \rangle\le m,
$$
or, in the codimension language,
$$
m\ge \codim\langle \bigcap\limits_E (Tf_\epsilon)_{u_E} \rangle
\ge \codim\langle (Tf_{\epsilon_0})_{u_E} \rangle
\ge \codim\langle \bigcup\limits_E (Tf_\epsilon)_{u_E} \rangle\ge 0.
$$

Maps $f_\delta\colon H_\delta\to M$, $\delta=1,\dots,\Delta$, are said to be
\begin{enumerate}
\item[\quad]{\it $\cap$-transversal}, if 
$$\codim\langle \bigcap\limits_E (Tf_\epsilon)_{u_E} \rangle=
\sum\limits_E \codim\langle (Tf_\epsilon)_{u_E} \rangle$$ for all
$E\subseteq\{1,\dots,\Delta\}$;
\item[\quad]{\it $\cup$-transversal}, if 
$$\sum\limits_E \dim\langle (Tf_\epsilon)_{u_E} \rangle=
\dim\langle \bigcup\limits_E (Tf_\epsilon)_{u_E} \rangle$$ for all
$E\subseteq\{1,\dots,\Delta\}$.
\end{enumerate}

The definition implies that $f_\delta$ can be $\cap$-transversal only for
$$\sum\limits_{\delta=1}^\Delta \codim\langle (Tf_\delta)_{u_{\{1,\dots,\Delta\}}} \rangle\le m$$
and, analogously,
$f_\delta$ can be $\cup$-transversal only for
$$\sum\limits_{\delta=1}^\Delta \dim\langle (Tf_\delta)_{u_{\{1,\dots,\Delta\}}} \rangle\le m.$$
It is easy to show that
$$
\sum\limits_{\delta=1}^\Delta \codim\langle (Tf_\delta)_{u_{\{1,\dots,\Delta\}}} \rangle\le m
\quad\textrm{and}\quad
\sum\limits_{\delta=1}^\Delta \dim\langle (Tf_\delta)_{u_{\{1,\dots,\Delta\}}} \rangle\le m
$$
comes into being simultaneously only for the case $\Delta=2$ and
$\codim\langle (Tf_1)_{u_{\{1,2\}}}\rangle+\codim\langle (Tf_2)_{u_{\{1,2\}}}\rangle=
\dim\langle (Tf_1)_{u_{\{1,2\}}}\rangle+\dim\langle (Tf_2)_{u_{\{1,2\}}}\rangle=m$. In this special case, concepts of
$\cap$-transversality and $\cup$-transversality are identical. (Sometimes, exactly this is understanded as a transversality:
cf. e.g. \cite{MAT}, Definition 7.5, where a decomposition of a tangent space onto two complementary subspaces is claimed.)

Let us consider $\cap$-transversal maps $f_\delta\colon H_\delta\to M$, $\delta=1,\dots,\Delta$ in the following situation:
$H_\delta$ are subsets (submanifolds) of $M$ and
$f_\delta\colon H_\delta\to M$ are their inclusion maps (immersions).
Then $H_\delta$ are called {\it $\cap$-transversal}, too. Moreover, if we have $\Delta$ foliations $F_\delta$ of $M$, we
take in every $u\in M$ their leaves: if they are {\it $\cap$-transversal} on each choice of $u$, we say that foliations $F_\delta$ of $M$ are {\it $\cap$-transversal}.

The concept {\it $\cup$-transversal} foliations $F_\delta$ of $M$ comes quite analogously.

A collection $\mathbf F=\{\mathcal F_\delta\}_{\delta=1}^\Delta$ of foliations of $M$ ($\dim M=m$)
with dimensions $p_\delta$ and codimensions $q_\delta$ 
is called the
{\it $\cap$-multifoliation} ({\it $\cup$-mul\-ti\-fo\-li\-a\-tion}), if
foliations $\mathcal F_\delta$ are 
$\cap$-transversal ($\cup$-transversal).
Especially, the $\cap$-multifoliation
($\cup$-multifoliation)
is called {\it total $\cap$-multifoliation}
({\it total $\cup$-multifoliation})
if $\Delta=m$.

It is clear that $q_1=\dots=q_\Delta=1$
for total $\cap$-multifoliation and $p_1=\dots=p_\Delta=1$
for total $\cup$-multifoliation.

\subsection
{Multifoliations by Kodaira and Spencer}

Kodaira and Spencer came in \cite{KOS} with the following concept of a multifoliation:
let $(P,\geq)$ be a partially ordered set with $\Delta$ elements and let us consider a surjective map 
$$
p\colon\{1,\dots,m\}\to P.
$$
(Thus, $m\ge\Delta$.)

We set 
$$
a^i_j=0
\;\;\textrm{for}\;\; 
p(j)\ngeq p(i)
$$ 
and denote by 
$$
\gel(m,\er;P,p)
$$
the subgroup of $\gel(m,\er)$ of linear transformations $\er^m\to\er^m$
given by $A=(a^i_j)$.
Further, we denote by $\Gamma(P,p)$ the pseudogroup of all local diffeomorphisms 
$$
g\colon U\to V,
\;\;
U,V\subseteq\er^m,
$$
such
that 
$$
\dej g_x\in \gel(m,\er;P,p)
\;\;\textrm{for all}\;\;
x\in U.
$$ 
A maximal atlas compatible with $\Gamma(P,p)$ is called the {\it $(P,p)$-multifoliate structure}
and $M$ endowed with a $(P,p)$-multifoliate structure is called the {\it $(P,p)$-multifoliate manifold}
or a manifold with a $(P,p)$-multifoliation.

\subsection
{Jets modulo multifoliations}

G. Ikegami has defined in his paper \cite{IKE} jets modulo foliations. 
(However, we refer also to approach of Doupovec, Kol\'a\v{r} and Mikulski, \cite{DOM}.)
We generalize Ikegami's concept by the following definition.
(In this section, we mean by a multifoliation either $\cap$-multifoliation or $\cup$-multifoliation or $(P,p)$-multifoliation.)

\begin{definition}
Let $H$, $M$ be two manifolds, $f,g\colon H\to M$ maps satisfying $f(h)=g(h)=u\in M$ and let $\mathbf F=\{\mathcal F_\delta\}_{\delta=1}^\Delta$ be a multifoliation of $M$. Then $f$ is said to have the {\it $(r_1,\dots,r_\Delta)$-multiorder contact modulo $\mathbf F$} with $g$ at $u$, if
for every $\Delta$-tuple of charts $\left\{U^\delta\ni u,\varphi^\delta\right\}_{1\le \delta\le \Delta}$ the maps
$$
\alpha^\delta\circ f\colon U^\delta\to\er^{p_\delta} \quad\textrm{and}\quad \alpha^\delta\circ g\colon U^\delta\to\er^{p_\delta}
$$
belong to the same (classical) $r_\delta$-jet at $u$.
(It means that for every curve $\gamma\colon\er\to H$ with $\gamma(0)=h$, the curves
$\alpha^\delta\circ f\circ\gamma$ and $\alpha^\delta\circ g\circ\gamma$ have the $r_\delta$-order contact in zero.)
As the relation "have the $(r_1,\dots,r_\Delta)$-multiorder contact modulo $\mathbf F$" is evidently an equivalence relation,
we denote the class of maps having the $(r_1,\dots,r_\Delta)$-multiorder contact modulo $\mathbf F$ with $f$ at $u$ by
$$
j^{r_1,\dots,r_\Delta}_h f \mathrm{mod}\mathbf F
$$
and call it {\it $(r_1,\dots,r_\Delta)$-jet modulo the multifoliation $\mathbf F$} with the {\it source $h\in H$} and the {\it target $u=f(h)\in M$}.
\end{definition}

We denote by $J^{r_1,\dots,r_\Delta}_h(H,M;\mathbf F)_u$ the set of all
$(r_1,\dots,r_\Delta)$-jets modulo the multifoliation $\mathbf F$ with the $h$ and the target $u$. Further, we denote
$$
J^{r_1,\dots,r_\Delta}_h(H,M;\mathbf F)=\bigcup\limits_{u\in M} J^{r_1,\dots,r_\Delta}_h(H,M;\mathbf F)_u,
$$
$$
J^{r_1,\dots,r_\Delta}(H,M;\mathbf F)_u=\bigcup\limits_{h\in H} J^{r_1,\dots,r_\Delta}_h(H,M;\mathbf F)_u
$$
and
$$
J^{r_1,\dots,r_\Delta}(H,M;\mathbf F)=\bigcup\limits_{u\in M}\bigcup\limits_{h\in H} J^{r_1,\dots,r_\Delta}_h(H,M;\mathbf F)_u.
$$

\begin{theorem}
For ma\-ni\-folds $H$ and $M$ and a mul\-ti\-foliation $\mathbf F$ of $M$,
spaces
$$
J^{r_1,\dots,r_\Delta}_h(H,M;\mathbf F)_u,
$$ 
$$
J^{r_1,\dots,r_\Delta}_h(H,M;\mathbf F),
$$
$$
J^{r_1,\dots,r_\Delta}(H,M;\mathbf F)_u
$$ 
and
$$J^{r_1,\dots,r_\Delta}(H,M;\mathbf F)$$
have a smooth manifold structure.
Further, there are bundle projections
$$
J^{r_1,\dots,r_\Delta}(H,M;\mathbf F)\to H \quad\textrm{and}\quad J^{r_1,\dots,r_\Delta}(H,M;\mathbf F)\to M
$$
as well as canonical bundle projections
$$
J^{r_1,\dots,r_\Delta}(H,M;\mathbf F)\to J^{\tilde r_1,\dots,\tilde r_\Delta}(H,M;\mathbf F)
$$
by restricting the multiorder, i.e. for
$0\le\tilde r_1\le r_1,\dots,0\le\tilde r_\Delta\le r_\Delta$. 
In doing so
$$
J^{0,\dots,0}(H,M;\mathbf F)=H\times M.
$$
\end{theorem}

Now, we present that $(R,S,Q)$-jet are included in 
the concept of the
$(r_1,\dots,r_\Delta)$-jet modulo the multifoliation $\mathbf F$. We recall
that two morphisms of fibered manifolds determine the same {\it $(R,S,Q)$-jet} ($R\le S$,
$R\le Q$) at a point $y$ if they have the same $R$-jet in $y$, their restrictions
to the fiber through $y$ have the same $S$-jet in $y$, and their base maps
have the same $Q$-jet in the base point of $y$. 

Let $Y\to M$ be a fibered manifold, $\dim M=q$, $\dim Y=p+q$.
The fibered manifold structure of $Y\to M$ determines
the foliation $\mathcal F_1$ with $p$-dimensional leaves (leaves = fibers).
Fiber bundles do not have global sections in general. However, if $Y\to M$ allow global sections (e.g. for a vector bundle with
smooth sections including zero section, such as constant smooth sections), we have also a foliation $\mathcal F_2$ with $q$-dimensional leaves (leaves = sections); of course, such $\mathcal F_2$ is not determined uniquely. But in the case we have obtained a (non-unique) multifoliation
which is simultaneously $\cup$-multifoliation and $\cap$-multifoliations.
Our construction implies:

\begin{theorem}
Let $\mathbf F$ be a multifoliation given by the fibration as stated above. Then there is a representation of every
$(R,S,Q)$-jet as a $(S,Q)$-jet modulo $\mathbf F$.
\end{theorem}

\section
{Final remarks}
\begin{remark}\rm{(On a unification.)}
Tom\'a\v{s} has generalized the concept of $(R,S,Q)$-jet to a concept of the nonholonomic $(R,S,Q)$-jet in \cite{TO1}.
He derived a composition of such jets as well as some properties of a corresponding bundle functor. That is the only attempt to unify
mentioned generalizations of jets up to now.
\end{remark}

\begin{remark}\rm{(Weil algebras approach.)}
Product preserving bundle functors on the category of smooth manifolds and smooth maps are in bijection with Weil algebras
and the natural transformations are in bijection with the Weil algebra homomorphisms.
The main example is the Weil algebra of the functor of $k$-dimensional $r$-th order velocities $J_0(\er^k,M)$.
For nonholonomic velocities or even quasivelocities, Weil algebra are described, see e.g. \cite{TOM}.
Further, it is known that product preserving bundle functors on the category 
fibered manifold are in bijection with algebra homorphisms acting between two Weil algebras.
These facts were recently also generalized. For $(P,p)$-multifoliate manifolds, Shurygin has defined in \cite{SHU}
an {\it inductive system of Weil algebra homomorphisms over $P$} as a collection $\mu=(A_\alpha,\mu^\beta_\alpha,P)$ consisting
of Weil algebras $A_\alpha$, $\alpha\in P$ and Weil algebra homomorphisms $\mu^\beta_\alpha\colon A_\beta\to A_\alpha$, $\beta\leq\alpha$
and has proved that product preserving bundle functors on $(P,p)$-multifoliate manifolds are uniquely determined by
inductive systems of Weil algebra homomorphisms.
\end{remark}

\begin{remark}\rm {(Weighted jets.)}
The concept of weighted jet bundles of sections of vector bundles over filtered manifolds was introduced by Morimoto (\cite{MOR})
in order to study differential equations on filtered manifold.
His weighted jet formalism is used in some problems connected with parabolic geometries (they represent special cases of Cartan's "espaces g\'en\'eralis\'es"
which are geometric structures that have homogeneous spaces $G/P$, where $G$ is a Lie group and $P$ a subgroup, as their models).

Kunzinger and Popovych (\cite{KUP}) identify differential consequences of a system of differential equations with a system of
algebraic equations in the jet space. For certail purposes, e.g. for different potential and pseudo-potential frames, they find useful to introduce the notion
of weight of differential variables instead of the order. Namely, for each variable in the infinite-order jet space
are weights of depending variables (coordinates in the target space) specified by the definition, with respect to a preservation of usual rules for
derivations of higher order jet coordinates.
\end{remark}

\subsection*{Acknowledgments}
Published results were acquired using the subsidization of the GA\v{C}R, grant No. 201/09/0981.
The literature recommendations of an unknown reviewer are acknowledged as well.

\end{document}